\documentclass[letterpaper, 10 pt, conference]{ieeeconf}  

\IEEEoverridecommandlockouts                              

\usepackage{graphicx}      
\usepackage{caption}
\usepackage{lineno,hyperref}
\usepackage{mathtools,amsfonts}
\usepackage{algorithm}
\usepackage{algpseudocode}
\usepackage{mathtools}
\usepackage{placeins}
\usepackage{color}
\usepackage[graphicx]{realboxes}
\usepackage{rotating}
\usepackage{afterpage}
\usepackage{subcaption}
\usepackage{lscape}
\usepackage{float}
\usepackage{lipsum}
\usepackage{rotating}
\usepackage{booktabs}
\usepackage{multirow}
\usepackage{siunitx}
\usepackage{svg}

\title{\LARGE \bf
Geometric Pairwise Classification of Affine Constraints over Ellipsoidal Regions*
}

\author{Igor M. L. Pataro$^{1}$, Martina Mammarella$^{2}$, Fabrizio Dabbene$^{2}$, and Teodoro Alamo$^{1}$
\thanks{*This publication is part of grant JDC2024-053231-I, funded by MICIU/AEI/10.13039/50110011033 and by the ESF+. T. Alamo acknowledges support from grant PID2025-173578OB-I00 funded by MICIU/AEI/10.13039/501100011033.}
\thanks{$^{1}$I.M.L. Pataro and T. Alamo are with Universidad de Sevilla, Escuela Superior de Ingenieros, Camino de los Descubrimientos s/n, Seville, Spain.
{\tt\small \{imendeslima,talamo\}@us.es.}}%
\thanks{$^{2}$M. Mammarella and F. Dabbene are with the CNR-IEIIT, c/o Politecnico di Torino, Corso Duca degli Abruzzi 24, Torino, Italy.
{\tt\small \{martina.mammarella,fabrizio.dabbene\}@cnr.it.}}%
\thanks{$\dagger$ Corresponding author: Igor M. L. Pataro {\tt\small imendeslima@us.es} }%
}

\begin{document}

\maketitle
\thispagestyle{empty}
\pagestyle{empty}

\begin{abstract}
This paper presents a deterministic geometric procedure for analysing affine constraints over fixed ellipsoidal regions, motivated by chance-constrained set approximation. Analytical tests classify inequalities as always violated, never violated, or cutting. Pairwise tests determine whether two inequalities can be strictly violated at a common point, yielding a compatibility matrix. The tests follow from quadratic optimisation and the Karush--Kuhn--Tucker conditions, with explicit treatment of tangencies, parallel normals, and degenerate inequalities. The formulation applies in arbitrary finite dimensions and provides geometric information for subsequent bounds on simultaneous sampled violations and safe-region approximation. 
\end{abstract}


\section{Introduction}
\label{sec:introduction}

Assessing constraint satisfaction under uncertainty is a recurring problem in optimisation, model calibration, and the validation of data-driven systems. Examples include constructing admissible parameter sets from uncertain measurements~\cite{Kosut1992,Mammarella2022} and evaluating the robustness of deep neural networks to random input perturbations~\cite{Weng2019PROVEN}. A common objective is to identify regions of the decision, parameter, or input space within which prescribed requirements hold with controlled risk. Chance constraints formalise this objective by limiting the probability of constraint violation to a prescribed level $\varepsilon\in(0,1)$~\cite{CharnesCooper1959}. This perspective motivates tractable descriptions of admissible regions that support both deterministic assessment and probabilistic certification. Whenever the requirements can be expressed or conservatively bounded by affine inequalities, their behaviour over an ellipsoidal region can be examined through geometric analysis, providing a building block for safe-region construction and chance-constrained set approximation.

Chance-constrained optimisation is generally challenging: evaluating violation probabilities may require multidimensional integration, and the resulting feasible set need not be convex~\cite{NemirovskiShapiro2006}. These difficulties motivate the construction of convex inner approximations that admit tractable deterministic representations. When the underlying chance-constrained set is time-invariant, such approximations can be computed offline and subsequently imposed as deterministic constraints in the online optimisation~\cite{Lorenzen2017}. Probabilistic scaling provides a systematic approach to this construction by separating the selection of an approximating set's shape and centre from the statistical calibration of its size~\cite{Alamo2019,Mammarella2022}. Under the corresponding sampling and calibration conditions, the resulting set is contained in the chance-constrained feasible region with a prescribed confidence. Ellipsoids are particularly convenient approximating sets because they provide a compact convex representation with an explicit description of their geometry.

In sampling-based formulations, fixing an uncertainty realisation can transform an uncertain constraint into a deterministic affine inequality of the form $c_i^\top x\leq d_i$. For a collection of such inequalities, assessing sampled violations over a candidate ellipsoid requires understanding how each constraint behaves throughout the region. An inequality may be satisfied everywhere, violated everywhere, or satisfied at some points and violated at others. In particular, an ellipsoid satisfying every sampled inequality is an inner approximation of the sampled feasible set. Establishing its relationship with the original chance-constrained set additionally requires statistical calibration. The geometric assessment is therefore a deterministic component of safe-region construction and violation analysis.

Motivated by this setting, this paper presents a geometric analysis of affine constraints over a fixed ellipsoidal region. Analytical expressions for the minimum and maximum of each affine function yield explicit conditions for classifying inequalities as always violated, never violated, or cutting. The resulting tests apply in arbitrary finite dimensions and avoid solving a separate numerical optimisation problem for each inequality. We also examine pairwise intersections of violation regions within the ellipsoid to construct a compatibility matrix, whose entries indicate whether two inequalities can be strictly violated at a common point. The same principle can be extended to groups of three, four, or up to $N$ inequalities by testing their joint violation regions, which requires information beyond pairwise compatibility.

The remainder of the paper is organised as follows. Section~\ref{sec:problem_formulation} states the problem and introduces the ellipsoidal region and affine inequalities. Section~\ref{sec:classification} develops the geometric classification and its analytical conditions. Section~\ref{sec:conclusions} presents the conclusions.

\section{Motivation and Problem Formulation}
\label{sec:problem_formulation}

Consider an optimisation problem in the decision vector $x\in\mathbb R^{n_x}$ subject to an affine chance constraint,
\begin{equation}
 \begin{array}{ll}
  \displaystyle\min_{x\in\mathcal X} & J(x)\\
  \text{subject to}
   & \operatorname{Viol}(x):=\Pr_w\{c(w)^\top x>d(w)\}\leq\varepsilon,
 \end{array}
 \label{eq:chance_problem}
\end{equation}
where $J$ is the objective function, $\mathcal X$ represents the deterministic admissible set, and $w$ is a random vector describing uncertainty. The measurable functions $c(w)\in\mathbb R^{n_x}$ and $d(w)\in\mathbb R$ define an inequality that is affine in $x$ for each uncertainty realisation. The parameter $\varepsilon\in(0,1)$ specifies the accepted violation risk. Denote the corresponding
chance-constrained feasible set by
\begin{equation}
 \mathcal F_\varepsilon
 =
 \{x\in\mathcal X:\operatorname{Viol}(x)\leq\varepsilon\}.
 \label{eq:chance_feasible_set}
\end{equation}

Following the set-approximation perspective in~\cite{Mammarella2022,Alamo2019}, we consider a candidate
ellipsoidal region 
\begin{equation}
 \mathcal E_\gamma
 =
 \left\{x\in\mathbb R^{n_x}:
 (x-x_c)^\top P^{-1}(x-x_c)\leq\gamma^2\right\},
 \label{eq:scaled_ellipsoid}
\end{equation}
where $x_c\in\mathbb R^{n_x}$ is its centre, $P=P^\top\succ0$ describes its shape, and $\gamma>0$ determines its scale. We assume
$\mathcal E_\gamma\subseteq\mathcal X$. The geometric analysis concerns a prescribed ellipsoid; for the statistical validation described below, its centre, shape, and scale are fixed independently of the calibration samples.

Let $w^{(1)},\ldots,w^{(N)}$ be independent and identically distributed samples drawn according to the distribution of $w$. Writing $c_i=c(w^{(i)})$ and $d_i=d(w^{(i)})$ gives the deterministic sampled inequalities
\begin{equation}
 c_i^\top x\leq d_i,
 \qquad i\in\mathcal I:=\{1,\ldots,N\}.
 \label{eq:sampled_affine_constraints}
\end{equation}
Their feasible half-spaces and common feasible region are
\begin{equation}
 \mathcal S_i=\{x:c_i^\top x\leq d_i\},
 \qquad
 \mathcal F_N=\mathcal X\cap
              \bigcap_{i\in\mathcal I}\mathcal S_i.
 \label{eq:sampled_feasible_set}
\end{equation}
Throughout this work, violation is strict: a point satisfying $c_i^\top x=d_i$ satisfies inequality $i$.

To relate these deterministic constraints to a statistical certificate, prescribe a confidence parameter $\delta\in(0,1)$ and an integer budget $0\leq r<N$. Select $N$ such that
\begin{equation}
 \textbf{B}(N,\varepsilon,r)
 :=
 \sum_{j=0}^{r}\binom Nj
       \varepsilon^j(1-\varepsilon)^{N-j}
 \leq\delta.
 \label{eq:binomial_calibration}
\end{equation}
Define the number of sampled inequalities that are violated somewhere within the ellipsoid as
\begin{equation}
 K_N(\mathcal E_\gamma)
 =
 \sum_{i=1}^{N}
 \mathbf 1_{\{\exists x\in\mathcal E_\gamma:
                    c_i^\top x>d_i\}}.
 \label{eq:set_violation_count}
\end{equation}
Applying the binomial validation argument in~\cite{Alamo2015} to the fixed candidate region gives
\begin{equation}
 \Pr_{\rm cal}
 \left\{
 K_N(\mathcal E_\gamma)\leq r,\quad
 \mathcal E_\gamma\not\subseteq\mathcal F_\varepsilon
 \right\}
 \leq\delta,
 \label{eq:ellipsoid_validation}
\end{equation}
where $\Pr_{\rm cal}$ denotes probability over the calibration samples. Indeed, the probability that a random inequality is violated somewhere in the ellipsoid bounds its violation probability at every individual point. Consequently, accepting a fixed ellipsoid when $K_N(\mathcal E_\gamma)\leq r$ provides a certificate with confidence at least $1-\delta$.

The deterministic task addressed in this paper is to determine
whether each sampled inequality holds throughout $\mathcal E_\gamma$, is violated everywhere within it, or is satisfied at some points and violated at others. These possibilities are characterised by the minimum and maximum of the affine residual $c_i^\top x-d_i$ over the ellipsoid, treating equality as constraint satisfaction. For pairs of inequalities, we additionally determine whether their strict violation regions intersect within $\mathcal E_\gamma$, meaning that a common point violates both constraints. The following section derives analytical geometric tests for individual classification and pairwise compatibility.

\section{CLASSIFICATION OF CONSTRAINTS}
\label{sec:classification}

We now consider a generic collection of sampled affine inequalities $c_i^\top x\leq d_i$, $i=1,\ldots,N$, and a reference ellipsoid constructed. Taking its centre and shape as given, we develop a deterministic geometric analysis applicable to the reference ellipsoid and its candidate enlargements. For any prescribed scale, the procedure classifies each inequality according to its violation region within the ellipsoid and identifies pairs that can be strictly violated at a common point. The analysis relies on the geometric relationship between the ellipsoid and each linear inequality, with the required distances and intersection points obtained through quadratic optimisation and the Karush--Kuhn--Tucker (KKT) conditions. The resulting pairwise compatibility information provides the basis for simultaneous constraint violation, which can be further extended to define deterministic bounds on maximum allowed violations for a scaled safe region.

\subsection{Geometry and joint distance computation}

Fix a scale $\gamma>0$ and consider the ellipsoid $\mathcal E_\gamma$ defined in \eqref{eq:scaled_ellipsoid}, with centre $x_c\in\mathbb R^{n_x}$ and shape matrix $P=P^\top\succ0$. Introducing the translated coordinates $z=x-x_c$ and the weighted norm $\|z\|_{P^{-1}}^2=z^\top P^{-1}z$, we equivalently write
\begin{equation}
 \mathcal E_\gamma
 =
 \left\{x_c+z:
 z\in\mathbb R^{n_x},
 \|z\|_{P^{-1}}\leq\gamma\right\}.
 \label{eq:cc_ellipsoid}
\end{equation}
The degenerate case $\gamma=0$ reduces to evaluating the inequalities at the single point $x_c$.

The sampled inequalities introduced in Section~\ref{sec:problem_formulation} are
\begin{equation}
 c_i^\top x\leq d_i,
 \qquad i\in\mathcal I:=\{1,\ldots,N\},
 \label{eq:cc_constraints}
\end{equation}
where $c_i\in\mathbb R^{n_x}$ and $d_i\in\mathbb R$. Their coefficients are treated as fixed throughout the following deterministic analysis. For $c_i\neq0$, the vector $c_i$ is the Euclidean normal to the boundary hyperplane. Let $C\in\mathbb R^{N\times n_x}$ have rows $c_i^\top$, and let $D=(d_1,\ldots,d_N)^\top$. Consistently with the strict-violation convention, define
\begin{equation}
 \mathcal V_i
 :=
 \{x\in\mathcal E_\gamma:c_i^\top x>d_i\},
 \label{eq:cc_violation_set}
\end{equation}
where dependence on the fixed scale $\gamma$ is suppressed for readability.

In translated coordinates, the boundary equation becomes
$c_i^\top z=\delta_i$, where
\begin{equation}
 \delta_i=d_i-c_i^\top x_c,
 \qquad
 \delta=D-Cx_c.
 \label{eq:cc_offsets}
\end{equation}
Partition the indices into
$\mathcal J_0:=\{i\in\mathcal I:c_i=0\}$ and
$\mathcal J_+:=\mathcal I\setminus\mathcal J_0$.
For all nonzero normals, the closest boundary points in the $P^{-1}$ metric can be computed jointly through the separable quadratic problem, for each $i\in\mathcal J_+$, as
\begin{equation}
 \begin{array}{ll}
 \displaystyle\min_{\{z_i\}_{i\in\mathcal J_+}}
 &
 \displaystyle\frac12
 \sum_{i\in\mathcal J_+}z_i^\top P^{-1}z_i
 \\[1mm]
 \mathrm{subject\ to}
 &
 c_i^\top z_i=\delta_i.
 \end{array}
 \label{eq:cc_joint_qp}
\end{equation}
Each variable $z_i\in\mathbb R^{n_x}$ corresponds to a separate hyperplane. Thus, the problem computes one projection per boundary; it does not require all boundary equations to hold at a common point. Each equality is feasible because $c_i\neq0$, and strict convexity ensures a unique minimiser. Zero-normal inequalities are classified separately in the following subsection.

The Lagrangian is
\begin{equation}
 \mathcal L
 =
 \sum_{i\in\mathcal J_+}
 \left[
 \frac12 z_i^\top P^{-1}z_i
 +\lambda_i(c_i^\top z_i-\delta_i)
 \right],
\end{equation}
where $\lambda_i\in\mathbb R$. Its KKT equations, for each $i\in\mathcal J_+$, separate into
\begin{equation}
 \begin{bmatrix}
 P^{-1}&c_i\\
 c_i^\top&0
 \end{bmatrix}
 \begin{bmatrix}
 z_i^\star\\
 \lambda_i^\star
 \end{bmatrix}
 =
 \begin{bmatrix}
 0\\
 \delta_i
 \end{bmatrix}.
 \label{eq:cc_single_kkt}
\end{equation}
Define
\begin{equation}
 g_i:=c_i^\top P c_i>0,
 \qquad q_i:=\sqrt{g_i}.
\end{equation}
Solving \eqref{eq:cc_single_kkt} gives
\begin{equation}
 \lambda_i^\star=-\frac{\delta_i}{g_i},
 \qquad
 z_i^\star=\frac{\delta_i}{g_i}Pc_i,
 \qquad
 x_i^\star=x_c+z_i^\star.
 \label{eq:cc_single_solution}
\end{equation}

From \eqref{eq:cc_single_solution}, the closest boundary point is $z_i^\star=(\delta_i/g_i)Pc_i$, with gives the squared distance from the centre
\begin{equation}
 \rho_i:=\|z_i^\star\|_{P^{-1}}^2
 =\frac{\delta_i^2}{g_i}.
 \label{eq:cc_single_distance}
\end{equation}

The boundary intersects the interior of $\mathcal E_\gamma$ when $\rho_i<\gamma^2$, is tangent when $\rho_i=\gamma^2$, and lies outside the ellipsoid when $\rho_i>\gamma^2$. Since $\rho_i$ is an unsigned squared distance, the sign of $\delta_i$ must also be retained to distinguish the feasible and violated sides.

\subsection{Exhaustive individual classification}

For the fixed scale $\gamma>0$ and each $i\in\mathcal J_+$, the extrema of the linear function $c_i^\top x$ over the ellipsoid are
\begin{equation}
 \begin{aligned}
 \min_{x\in\mathcal E_\gamma}c_i^\top x
 &=c_i^\top x_c-\gamma q_i,\\
 \max_{x\in\mathcal E_\gamma}c_i^\top x
 &=c_i^\top x_c+\gamma q_i.
 \end{aligned}
 \label{eq:cc_extrema}
\end{equation}
Indeed, Cauchy--Schwarz gives
\[
 |c_i^\top(x-x_c)|
 \leq q_i\|x-x_c\|_{P^{-1}}
 \leq\gamma q_i,
\]
with the minimum and maximum attained at
$x_c-\gamma Pc_i/q_i$ and
$x_c+\gamma Pc_i/q_i$, respectively.
Using $\delta_i=d_i-c_i^\top x_c$ and
$\rho_i=\delta_i^2/q_i^2$, the nonzero-normal inequalities are classified as
\begin{equation}
 \begin{array}{ll}
 \mathrm{always\ violated:}
 &\delta_i<0,\quad\rho_i>\gamma^2,\\
 \mathrm{never\ violated:}
 &\delta_i\geq0,\quad\rho_i\geq\gamma^2,\\
 \mathrm{cutting:}
 &\rho_i<\gamma^2\ \mathrm{or} \:\:
   (\delta_i<0,\ \rho_i=\gamma^2).
 \end{array}
 \label{eq:cc_nonzero_classes}
\end{equation}
Equivalently, these three cases correspond to
$\delta_i<-\gamma q_i$,
$\delta_i\geq\gamma q_i$, and
$-\gamma q_i\leq\delta_i<\gamma q_i$,
respectively.

When $\rho_i<\gamma^2$, the boundary hyperplane intersects the interior of $\mathcal E_\gamma$. At $\delta_i=\gamma q_i$, the boundary is tangent to the ellipsoid and no point strictly violates the inequality. At $\delta_i=-\gamma q_i$, every point except the unique tangency point strictly violates it. This latter case is included in the cutting class because its violation region is nonempty but does not contain the entire ellipsoid.

For $i\in\mathcal J_0$, the inequality reduces to $0\leq d_i$. Consequently, $d_i<0$ implies $\mathcal V_i=\mathcal E_\gamma$, whereas $d_i\geq0$ implies $\mathcal V_i=\emptyset$. These inequalities are classified directly, without computing a distance or solving a boundary KKT system.

Including both zero and nonzero normals, define
\begin{equation}
 \begin{aligned}
 \mathcal I_{\rm A}
 &:=\{i\in\mathcal I:\mathcal V_i=\mathcal E_\gamma\},\\
 \mathcal I_{\rm N}
 &:=\{i\in\mathcal I:\mathcal V_i=\emptyset\},\\
 \mathcal I_{\rm C}
 &:=\mathcal I\setminus
       (\mathcal I_{\rm A}\cup\mathcal I_{\rm N}).
 \end{aligned}
 \label{eq:cc_index_classes}
\end{equation}
These sets form a disjoint partition of $\mathcal I$. For every $i\in\mathcal I_{\rm C}$, the normal is nonzero, the violation region $\mathcal V_i$ is nonempty, and $\rho_i\leq\gamma^2$. Hence, the projected boundary point $x_i^\star$ defined in \eqref{eq:cc_single_solution} belongs to $\mathcal E_\gamma$, although it does not itself strictly violate inequality $i$. Only the indices in $\mathcal I_{\rm C}$ require pairwise compatibility analysis; the always-violated inequalities contribute the fixed count $|\mathcal I_{\rm A}|$, while the never-violated inequalities contribute zero.

\subsection{Compatibility matrix and direct decisions}

Let $m=|\mathcal I_{\rm C}|$ and define a symmetric Boolean matrix $M\in\{0,1\}^{m\times m}$. Fix an ordering of $\mathcal I_{\rm C}$ and label the rows and columns of $M$ by the corresponding original constraint indices. For $i,j\in\mathcal I_{\rm C}$, define
\begin{equation}
\begin{array}{ll}
 M_{ii}=0,\\
 M_{ij}=1
 &\ \Longleftrightarrow
 \mathcal V_i\cap\mathcal V_j\neq\emptyset,
 \qquad i\neq j.
 \end{array}
 \label{eq:cc_M_definition}
\end{equation}
Thus, an off-diagonal entry equals one precisely when a point in $\mathcal E_\gamma$ strictly violates both inequalities. An intersection of their boundary hyperplanes alone is insufficient.

Only cutting inequalities are represented in $M$. Set all entries of $M$ to zero, then process distinct pairs in $\mathcal I_{\rm C}$. If $\mathcal I_{\rm C}$ is empty, $M$ is the empty matrix and no pairwise analysis is required.
\hfill \break
\hfill \break
\noindent\textit{Common centre witness:}
Since $x_c\in\mathcal E_\gamma$ and
$c_i^\top x_c>d_i\Longleftrightarrow\delta_i<0$, define
\begin{equation}
 \mathcal I_c
 :=
 \{i\in\mathcal I_{\rm C}:\delta_i<0\}.
 \label{eq:cc_center_group}
\end{equation}
The centre strictly violates every inequality indexed by $\mathcal I_c$. Therefore, all distinct pairs within this set can be marked simultaneously:
$M_{ij}=M_{ji}=1$ for $i,j\in\mathcal I_c$, $i\neq j$.
Pairs with only one index in $\mathcal I_c$ remain unresolved by this test.
\hfill \break
\hfill \break
\noindent\textit{Parallel normals:}
Consider an unresolved pair $i,j\in\mathcal I_{\rm C}$ with
$c_j=\alpha c_i$, where $\alpha\neq0$.
Both normals are nonzero, and
$q_j=|\alpha|q_i$.
Introducing the scalar
\[
 s=c_i^\top(x-x_c),
\]
the extrema in \eqref{eq:cc_extrema} imply that $s$ attains every value in
$[-\gamma q_i,\gamma q_i]$ over $\mathcal E_\gamma$.
Simultaneous strict violation is therefore equivalent to finding an admissible $s$ satisfying
\[
 s>\delta_i,
 \qquad
 \alpha s>\delta_j.
\]

If $\alpha>0$, the normals have the same orientation, and both violations require
\begin{equation}
 s>\max\{\delta_i,\delta_j/\alpha\}.
 \label{eq:cc_same_direction}
\end{equation}
Because both inequalities are cutting,
$\delta_i<\gamma q_i$ and
$\delta_j<\gamma q_j=\gamma\alpha q_i$.
Hence both thresholds are strictly smaller than $\gamma q_i$, so an admissible $s$ exists and $M_{ij}=M_{ji}=1$.
This includes coincident boundary hyperplanes with the same orientation.

If $\alpha<0$, write $\beta=-\alpha>0$.
The violations instead require
\[
 \delta_i<s<-\delta_j/\beta.
\]

For two cutting inequalities,
$\delta_i\geq-\gamma q_i$ and
$\delta_j\geq-\gamma q_j=-\gamma\beta q_i$.
Consequently, the condition to define compatibility between the pairs is resumed as
\[
 M_{ij}=1
 \quad\Longleftrightarrow\quad
 \delta_i<-\delta_j/\beta.
\]
In particular, coincident boundary hyperplanes with opposite orientations satisfy $\delta_j=-\beta\delta_i$ and cannot be strictly violated simultaneously.
\hfill \break
\hfill \break

\subsection{Independent normals: point tests and intersection distance}

Consider an unresolved pair $i,j\in\mathcal I_{\rm C}$ with linearly independent normals. The boundary projections $x_i^\star$ and $x_j^\star$ defined in \eqref{eq:cc_single_solution} belong to $\mathcal E_\gamma$. First evaluate
\begin{equation}
 c_j^\top x_i^\star>d_j
 \quad\mathrm{or}\quad
 c_i^\top x_j^\star>d_i.
 \label{eq:cc_point_tests}
\end{equation}
Either condition implies $M_{ij}=M_{ji}=1$. To justify the first condition, note that $x_i^\star$ lies on boundary $i$ and strictly violates inequality $j$. 

If neither point test succeeds, define
\begin{equation}
 B_{ij}
 :=
 \begin{bmatrix}
 c_i^\top\\
 c_j^\top
 \end{bmatrix},
 \qquad
 \Delta_{ij}
 :=
 \begin{bmatrix}
 \delta_i\\
 \delta_j
 \end{bmatrix},
 \label{eq:cc_pair_data}
\end{equation}
where $B_{ij}\in\mathbb R^{2\times n_x}$ and
$\Delta_{ij}\in\mathbb R^2$. The closest point to the intersection of the two boundary hyperplanes, measured from $x_c$ in the $P^{-1}$ metric, is obtained through
\begin{equation}
 \begin{array}{ll}
 \displaystyle\min_{z\in\mathbb R^{n_x}}
 &\displaystyle\frac12 z^\top P^{-1}z\\[1mm]
 \mathrm{subject\ to}
 &B_{ij}z=\Delta_{ij}.
 \end{array}
 \label{eq:cc_pair_qp}
\end{equation}
Linear independence of the normals gives $B_{ij}$ full row rank, ensuring feasibility of the equalities. Together with $P\succ0$, this also ensures that the KKT system
\begin{equation}
 \begin{bmatrix}
 P^{-1}&B_{ij}^\top\\
 B_{ij}&0_{2\times2}
 \end{bmatrix}
 \begin{bmatrix}
 z_{ij}^\star\\
 \lambda_{ij}^\star
 \end{bmatrix}
 =
 \begin{bmatrix}
 0\\
 \Delta_{ij}
 \end{bmatrix}
 \label{eq:cc_pair_kkt}
\end{equation}
is nonsingular, where $\lambda_{ij}^\star\in\mathbb R^2$.

Define the positive definite Gram matrix
\[
 G_{ij}
 :=
 B_{ij}PB_{ij}^\top.
\]
The analytical solution is
\begin{equation}
 \begin{aligned}
 z_{ij}^\star
 &=PB_{ij}^\top G_{ij}^{-1}\Delta_{ij},\\
 \lambda_{ij}^\star
 &=-G_{ij}^{-1}\Delta_{ij},\\
 x_{ij}^\star
 &=x_c+z_{ij}^\star.
 \end{aligned}
 \label{eq:cc_pair_solution}
\end{equation}
For $n_x=2$, $x_{ij}^\star$ is the unique intersection of the two boundary lines. For $n_x>2$, it is the closest point of their affine intersection in the $P^{-1}$ metric. Its squared distance from the centre is
\begin{equation}
 \begin{aligned}
 \rho_{ij}
 &:=(z_{ij}^\star)^\top P^{-1}z_{ij}^\star\\
 &=\Delta_{ij}^\top G_{ij}^{-1}\Delta_{ij}.
 \end{aligned}
 \label{eq:cc_pair_distance}
\end{equation}
For an explicit scalar evaluation, set
$g_{ij}:=c_i^\top Pc_j$. Then
\[
 G_{ij}
 =
 \begin{bmatrix}
 g_i&g_{ij}\\
 g_{ij}&g_j
 \end{bmatrix},
 \qquad
 g_i g_j-g_{ij}^2>0,
\]
and
\begin{equation}
 \rho_{ij}
 =
 \frac{
 g_j\delta_i^2
 -2g_{ij}\delta_i\delta_j
 +g_i\delta_j^2
 }{
 g_i g_j-g_{ij}^2
 }.
 \label{eq:cc_pair_scalar}
\end{equation}

\begin{figure*}
    \centering
    \includegraphics[width=1\linewidth,trim={10mm 10mm 10mm 10mm}]{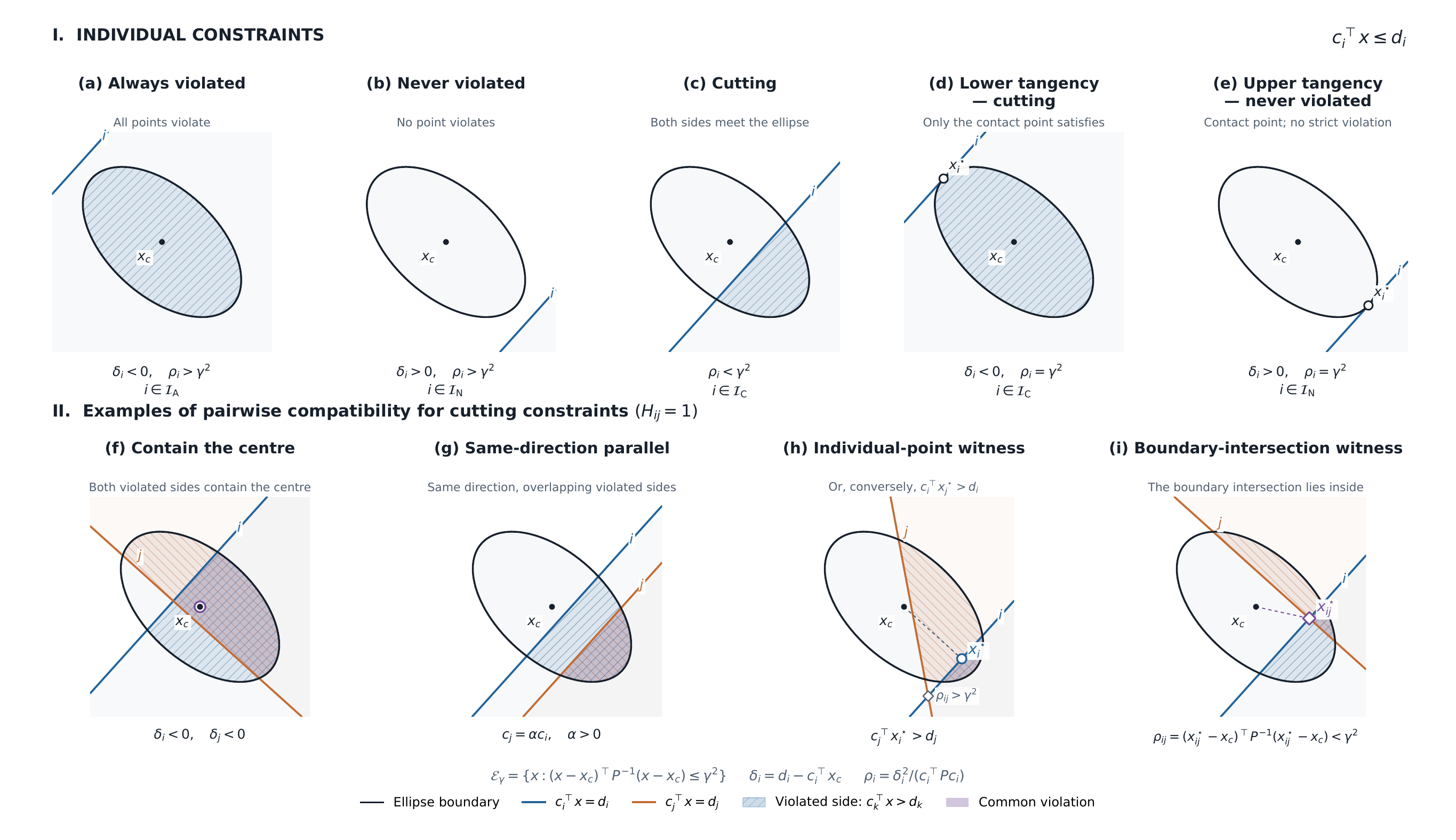}
    \caption{Constraint classification for the ellipsoid safe region.}
    \label{fig:constraint_clasf}
\end{figure*}

\noindent\textbf{Proposition 1:}

For a pair $i,j\in\mathcal I_{\rm C}$ with linearly independent normals, if neither the common-centre test nor the individual point tests identify compatibility, then
\begin{equation}
 M_{ij}=1
 \quad\Longleftrightarrow\quad
 \rho_{ij}<\gamma^2.
 \label{eq:cc_final_pair_test}
\end{equation}

\noindent\textit{Proof:}
Consider the closed convex set
\[
 \mathcal Q_{ij}
 :=
 \{z\in\mathbb R^{n_x}:B_{ij}z\geq\Delta_{ij}\},
\]
where the inequalities are componentwise. Full row rank of $B_{ij}$ ensures that $\mathcal Q_{ij}$ is nonempty. Let $\widehat z_{ij}$ be its unique minimum-norm point in the $P^{-1}$ metric.

Suppose $M_{ij}=1$. A point strictly violating both inequalities exists in $\mathcal E_\gamma$. If necessary, moving it slightly towards $x_c$ preserves both strict violations and places it in the ellipsoid interior. Consequently,
\[
 \|\widehat z_{ij}\|_{P^{-1}}^2<\gamma^2.
\]
If neither inequality is active at $\widehat z_{ij}$, stationarity gives $\widehat z_{ij}=0$, so the common-centre test succeeds. If exactly one inequality is active, the KKT conditions give $\widehat z_{ij}=z_i^\star$ or $\widehat z_{ij}=z_j^\star$, while the other inequality is strict. The corresponding individual point test therefore succeeds. Both cases are excluded by hypothesis. Hence, both inequalities are active, and the KKT conditions coincide with \eqref{eq:cc_pair_kkt}. Thus,
$\widehat z_{ij}=z_{ij}^\star$ and $\rho_{ij}<\gamma^2$.

Conversely, suppose $\rho_{ij}<\gamma^2$. Since $B_{ij}$ has full row rank, there exists a direction $v$ satisfying
$B_{ij}v=(1,1)^\top$. For sufficiently small $t>0$, the point
\[
 x=x_c+z_{ij}^\star+tv
\]
remains in $\mathcal E_\gamma$ and satisfies
\[
 B_{ij}(x-x_c)
 =
 \Delta_{ij}+t(1,1)^\top
 >
 \Delta_{ij}.
\]
It therefore strictly violates both inequalities, proving $M_{ij}=1$.

Under the preceding exclusions, $\rho_{ij}=\gamma^2$ represents contact without simultaneous strict violation, whereas $\rho_{ij}>\gamma^2$ implies incompatibility. The condition $\rho_{ij}<\gamma^2$ is sufficient for compatibility for any independent pair; its necessity requires the preceding tests. In particular, simultaneous violations may exist even when the intersection of the two boundary hyperplanes lies outside the ellipsoid.

\subsection{Assembly and implication for violation counting}

The procedure to count the simultaneous violation first computes the offsets $\delta$, evaluates the analytical solution of \eqref{eq:cc_single_kkt}, and partitions the inequalities into $\mathcal I_{\rm A}$, $\mathcal I_{\rm N}$, and $\mathcal I_{\rm C}$. The geometric cases considered throughout this section are illustrated in Fig.~\ref{fig:constraint_clasf}. Always-violated inequalities are counted separately, while never-violated inequalities require no further analysis. The compatibility matrix $M$ is constructed exclusively from the cutting inequalities indexed by $\mathcal I_{\rm C}$.

After initialising $M$ to zero, all distinct pairs within the common-centre group $\mathcal I_c$ are marked as compatible. The remaining unordered pairs $i,j\in\mathcal I_{\rm C}$, $i<j$, are examined using the parallel-normal conditions. For independent normals, the individual point tests are applied first. If neither succeeds, compatibility is determined using the analytical pairwise squared distance and the condition $\rho_{ij}<\gamma^2$. Each decision is stored symmetrically as $M_{ij}=M_{ji}$, and the diagonal remains zero. At most $|\mathcal I_{\rm C}|(|\mathcal I_{\rm C}|-1)/2$ unordered pairs are considered, of which $|\mathcal I_c|(|\mathcal I_c|-1)/2$ are resolved by the common-centre test. The closed-form expressions in \eqref{eq:cc_pair_kkt} avoid solving the full optimisation problem in \eqref{eq:cc_pair_qp} by using the KKT conditions, making the problem computationally tractable.

The resulting matrix records which pairs of cutting inequalities can be strictly violated at a common point in $\mathcal E_\gamma$. Every always-violated inequality contributes one violation at every point, whereas never-violated inequalities contribute zero.
Consequently,
\begin{equation}
 \Gamma_N(\gamma)
 =
 |\mathcal I_{\rm A}|
 +
 \max_{x\in\mathcal E_\gamma}
 \sum_{i\in\mathcal I_{\rm C}}
 \mathbf 1_{\{c_i^\top x>d_i\}}.
 \label{eq:cc_count_decomposition}
\end{equation}

This decomposition and the compatibility matrix $M$ provide a geometric basis for analysing simultaneous sampled violations and developing computable deterministic upper bounds on $\Gamma_N(\gamma)$.

The preceding tests assume exact arithmetic. In numerical implementations, floating-point errors can affect comparisons near the classification and compatibility thresholds. These cases may be handled conservatively by introducing a small comparison tolerance in  \eqref{eq:cc_nonzero_classes} and
\eqref{eq:cc_final_pair_test}. For example, after normalisation, a dimensionless tolerance $\tau=10^{-8}$ may be used to replace the strict pairwise test $\rho_{ij}/\gamma^2<1$ by
\[ \rho_{ij}/\gamma^2\leq1+\tau. \]

Similarly, inequalities whose classification is uncertain near a tangency threshold are retained as candidates rather than discarded. This convention may retain additional inequalities or mark additional pairs as compatible, making the resulting violation bound more conservative. The tolerance must be chosen based on numerical scaling and conditioning.

Algorithm~\ref{alg:constraint_classification} first partitions the sampled inequalities and counts those with nonempty violation regions, giving $K_N(\mathcal E_\gamma)=n_{\rm A}+n_{\rm C}$. The threshold conditions also classify zero-normal inequalities directly; boundary projections are computed only for cutting inequalities, for which $g_i>0$. The pairwise stage constructs $M$ exclusively from these cutting inequalities, preserving symmetry and a zero diagonal. Always-violated inequalities contribute $n_{\rm A}$ at every point, while the compatibility matrix provides the geometric information required for subsequent bounds on simultaneous violations.

\begin{algorithm}[!t]
\caption{Constraint counting and pairwise classification}
\label{alg:constraint_classification}
\small
\begin{algorithmic}[1]
\Require Inequalities $\{(c_i,d_i)\}_{i=1}^{N}$;
         centre $x_c$; $P=P^\top\succ0$; scale $\gamma>0$.
\Ensure Classes $\mathcal I_{\rm A},\mathcal I_{\rm N},
         \mathcal I_{\rm C}$, their counts, and compatibility
         matrix $M$.

\Statex \textit{Individual classification}
\State $\mathcal I_{\rm A}\gets\emptyset$,
       $\mathcal I_{\rm N}\gets\emptyset$,
       $\mathcal I_{\rm C}\gets\emptyset$
\For{$i=1,\ldots,N$}
    \State $\delta_i\gets d_i-c_i^\top x_c$
    \State $p_i\gets Pc_i$,
           $g_i\gets c_i^\top p_i$,
           $q_i\gets\sqrt{g_i}$
    \If{$\delta_i<-\gamma q_i$}
        \State $\mathcal I_{\rm A}
                \gets\mathcal I_{\rm A}\cup\{i\}$
    \ElsIf{$\delta_i\geq\gamma q_i$}
        \State $\mathcal I_{\rm N}
                \gets\mathcal I_{\rm N}\cup\{i\}$
    \Else
        \State $\mathcal I_{\rm C}
                \gets\mathcal I_{\rm C}\cup\{i\}$
        \State $x_i^\star\gets
                x_c+(\delta_i/g_i)p_i$
    \EndIf
\EndFor
\State $n_{\rm A}\gets|\mathcal I_{\rm A}|$,
       $n_{\rm N}\gets|\mathcal I_{\rm N}|$,
       $n_{\rm C}\gets|\mathcal I_{\rm C}|$
\State $K_N(\mathcal E_\gamma)\gets n_{\rm A}+n_{\rm C}$

\Statex \textit{Initialisation and common-centre test}
\State $M\gets0_{n_{\rm C}\times n_{\rm C}}$,
       with rows and columns labelled by $\mathcal I_{\rm C}$
\State $\mathcal I_c\gets
       \{i\in\mathcal I_{\rm C}:\delta_i<0\}$
\State Set $M_{ij}=M_{ji}=1$ for all
       $i,j\in\mathcal I_c$ with $i<j$

\Statex \textit{Remaining pairwise decisions}
\ForAll{$i,j\in\mathcal I_{\rm C}$ with
         $i<j$ and $M_{ij}=0$}
    \If{$c_j=\alpha c_i$ for some $\alpha\neq0$}
        \If{$\alpha>0$}
            \State $m\gets1$
        \Else
            \State $m\gets
             \mathbf 1_{\{\delta_i<\delta_j/\alpha\}}$
        \EndIf
    \ElsIf{$c_j^\top x_i^\star>d_j$
            or $c_i^\top x_j^\star>d_i$}
        \State $m\gets1$
    \Else
        \State $B_{ij}\gets
               \begin{bmatrix}c_i&c_j\end{bmatrix}^{\top}$,
               $\Delta_{ij}\gets
               \begin{bmatrix}\delta_i&\delta_j\end{bmatrix}^{\top}$
        \State $G_{ij}\gets B_{ij}PB_{ij}^{\top}$
        \State $\rho_{ij}\gets
               \Delta_{ij}^{\top}G_{ij}^{-1}\Delta_{ij}$
        \State $m\gets
               \mathbf 1_{\{\rho_{ij}<\gamma^2\}}$
    \EndIf
    \State $M_{ij}\gets h$, $M_{ji}\gets h$
\EndFor
\State \Return Classes, their counts,
       $K_N(\mathcal E_\gamma)$, and $M$
\end{algorithmic}
\end{algorithm}

\section{CONCLUSIONS}
\label{sec:conclusions}

This paper presented a deterministic geometric procedure for classifying affine constraints over ellipsoidal regions. Analytical tests distinguish always-violated, never-violated, and cutting inequalities, while pairwise analysis determines whether two constraints can be strictly violated at a common point. The resulting compatibility matrix and individual counts provide information for bounding simultaneous sampled violations and supporting probabilistically calibrated set approximations. The formulation applies in arbitrary finite dimensions and accounts consistently for boundary contact and degenerate cases. Future work will combine the proposed geometric classification and pairwise violation tests with graph-based methods, including graph colouring, $m$-core decompositions, and maximum-clique algorithms, to obtain computable bounds on the maximum number of constraints that can be strictly violated simultaneously at a single point within the ellipsoid.

\addtolength{\textheight}{-12cm}   







\bibliographystyle{IEEEtran}
\bibliography{references.bib}

\end{document}